\documentclass[12pt]{amsart}

\usepackage{amsthm}
\usepackage{hyperref}
\usepackage{geometry}
\usepackage{xpatch}
\usepackage{amssymb}
\usepackage{setspace}
\usepackage{amsthm}
\usepackage{enumerate}
\usepackage{commath}
\usepackage{hyperref}
\usepackage{titlesec}
\usepackage{amsmath,amssymb,amsfonts,anysize,mathrsfs,graphicx,setspace,parskip,multicol,upgreek,yfonts,tikz,physics,graphicx,xcolor,enumerate,csquotes,hyperref}

\newtheoremstyle{mystyle}
  {}
  {}
  {\itshape}
  {}
  {\bfseries}
  {}
  {.5em}
  {}

\theoremstyle{mystyle}

\newtheorem{theorem}{Theorem}
\newtheorem{problem}{Problem}

\newtheorem{claim}{Claim}

\newtheorem{lemma}{Lemma}

\newtheorem{corollary}{Corollary}
\titleformat{\section}{\small\bfseries}{\thesection}{1em}{}

\theoremstyle{definition}
\newtheorem{definition}{Definition}
\newtheorem{example}{Example}

\numberwithin{equation}{section}
\begin{document}

\title{Semi-proximal Spaces and Normality}
\author[K. Almontashery]{Khulod Almontashery} 
\address{Department of Mathematics and Statistics\\ York University\\ Toronto, ON M3J 1P3\\Canada}
\email{khulod@yorku.ca}

\author[P.J. Szeptycki]{Paul J. Szeptycki}
\address{Department of Mathematics and Statistics\\ York University\\ Toronto, ON M3J 1P3\\Canada}
\email{szeptyck@yorku.ca}

\begin{abstract}
We consider the relationship between normality and semi-proximality. We give a consistent example of a first countable locally compact Dowker space that is not semi-proximal, and two ZFC examples of semi-proximal non-normal spaces. This answers a question of Nyikos. One of the examples is a subspace of $(\omega+1) \times \omega_1$. In contrast, we show that every normal subspace of a finite power of $\omega_1$ is semi-proximal. 
\end{abstract}



\maketitle
\section{Introduction}

Suppose $(X, {\mathfrak U})$ is a uniform space, where ${\mathfrak U}$ is either a uniformity or a uniformity base. The proximal game is defined on $(X,{\mathfrak U})$ as follows. In inning 0, Player I chooses an entourage $U_0$ and Player II chooses $x_0 \in X$. In inning $n + 1$, Player I chooses an entourage $U_{n+1} \subseteq U_n$ and Player II chooses $x_{n+1} \in U_n[x_n] = \{y\in X : \langle x_n,y\rangle\in U_n\}$. Player I wins the game if either  $\bigcap_{n<\omega} U_n[x_n] = \emptyset$ or the sequence $(x_n : n <\omega)$ converges.

\begin{definition} A space $X$ is {\it proximal} if there is a compatible uniformity ${\mathfrak U}$ on $X$ such that there is a winning strategy for Player I in the proximal game on $(X,{\mathfrak U})$. \end{definition}

The proximal game and the class of proximal spaces were originally introduced by J. Bell in \cite{B1} as a means of studying normality in uniform box products. 

Naturally, in any topological game, one Player having no winning strategy is a weaker notion than the other Player having a winning strategy, and Nyikos defined the corresponding class for the proximal game as the class of semi-proximal spaces in \cite{Ny}:

\begin{definition}  A space $X$ is {\it semi-proximal} if there is a compatible uniformity ${\mathfrak U}$ on $X$ such that there is no winning strategy for Player II in the proximal game on $(X,{\mathfrak U})$.  \end{definition}

And so, evidently, any proximal space is semi-proximal. In addition to introducing this class, Nyikos proved, among other results, that semi-proximal spaces are Fr\'echet.

Nyikos posed the question of whether semi-proximal spaces are normal (Problem 13 of \cite{Ny}). We will present two counter-examples to Nyikos's question and consider some other closely related questions concerning more broadly the relationship between the class of normal spaces and the class of semi-proximal spaces. 

For example, if $X$ is semi-proximal and $Y$ is proximal, then the product $X\times Y$ is also semi-proximal \cite{Pa}. So in particular, $X\times (\omega+1)$ is then semi-proximal. Recall that a Dowker space is a normal space whose product with $\omega+1$ is not normal. So a counter-example to Nyikos's question could also be obtained by constructing a semi-proximal Dowker space. On the other hand, since semi-proximal spaces are Fr\'echet we can deduce that not every normal space is semi-proximal. For example, consider $\omega_1+1$ or any other normal non-Fr\'echet space.  However, for the class of Fr\'echet spaces the question is more interesting. While we were not able to construct a semi-proximal Dowker space, we do have a consistent example of first countable, locally compact (hence Fr\'echet and much more) Dowker space that is not semi-proximal. We do not know of any other first countable or even Fr\'echet normal space that is not semi-proximal. 

Characterizing normality in subspaces of products of ordinals has been studied extensively (e.g., \cite{KS}, \cite{Fl}, \cite{S}). For example, Fleissner proved that a subspace of a finite product of ordinals is normal if and only if  it has property ${\mathcal P}$, where ${\mathcal P}$ is either collectionwise normal, normal and strongly zero-dimensional, or shrinking.

In the last section we show that normality and semi-proximality are equivalent in finite product of subspaces of $\omega_1$. Also we prove that normal subspaces of finite powers of $\omega_1$ are semi-proximal, but the converse fails as we construct an example of a semi-proximal non-normal subspace of $(\omega+1) \times \omega_1$. 

We refer the reader to \cite{Eng} for the terminology and basic theory of uniformities and uniform spaces. Recall that uniformities have an equivalent formulation in terms of normal families of open covers. So if ${\mathcal A}$ is the corresponding normal family of open covers corresponding to the uniformity ${\mathfrak U}$ then the proximal game is equivalently described as follows: in inning 0 Player I chooses an open cover $A_0\in {\mathcal A}$ and Player II chooses $x_0 \in X$. In inning $n + 1$, Player I chooses a ${A}_{n+1}\in {\mathcal A}$ that refines ${A}_n$ and Player II chooses $x_{n+1} \in \,\,\text{St}(x_n,{ A}_n) = \bigcup\{U\in { A}_n :  x_n\in U\}$. Player I wins the game if either  $\bigcap_{n<\omega}\,\, \text{St}(x_n,{ A}_n)  = \emptyset$ or the sequence $(x_n : n <\omega)$ converges. 

It is easy to check that if $X = \bigoplus_{i\in I} X_i$, and $X_i$ is semi-proximal for all $i\in I$, then $X$ is semi-proximal. Also being proximal or semi-proximal are both closed-hereditary (e.g., see \cite{B1}).

Our notation and terminology are standard -- see \cite{Eng} for any topological notions and \cite{K} for set-theoretic notions. For background on Dowker spaces, and in particular the construction of deCaux-type examples see \cite{RM}.

\section{A not semi-proximal Dowker space}

We describe a de Caux type Dowker space constructed by enhancing $\clubsuit$ to a stronger principle we denote $\clubsuit^\ast$. 

\begin{definition} The principle $\clubsuit^\ast$ is the statement that there exists a sequence $\{C_\alpha : \alpha\in\text{Lim}(\omega_1)\}$, where $C_\alpha\subset \alpha$, has order type $\omega$, and $\sup(C_\alpha) = \alpha$ such that for every uncountable subset $A$ of $\omega_1$, $\{\alpha<\omega_1 : |A\cap C_\alpha| = \omega\}$ contains a club.
\end{definition}

It is easy to see that $\clubsuit^\ast$ follows from $\diamondsuit^\ast$ (see \cite{K}). Indeed, if $\{{\mathcal A}_\alpha:\alpha<\omega_1\}$ is a $\diamondsuit^\ast$ sequence and if for each limit $\alpha$, $C_\alpha$ is an increasing sequence cofinal with $\alpha$  satisfying $C_\alpha \cap A$ is infinite for all $A\in {\mathcal A}_\alpha$, then $\{C_\alpha:\alpha<\omega_1\}$ will be a $\clubsuit^\ast$ sequence.

 Our space is just one of the standard de Caux spaces constructed from our $\clubsuit^*$ sequence $\{C_\alpha:\alpha<\omega_1\}$: it has as its underlying set $X = \omega_1\times \omega$, and the topology is defined so that 
\begin{enumerate}
\item $C_\alpha\times\{n-1\}$ converges to $\langle\alpha,n\rangle$ for each limit ordinal $\alpha<\omega_1$. 
\item $\langle\alpha,n\rangle$ is isolated for all successor ordinals $\alpha$. 
\item $\alpha+1\times \omega$ is clopen for all $\alpha\in \omega_1$
\item The space is first countable and locally compact.
\end{enumerate}

One way to describe the topology is to just declare the family $\{\{\langle\alpha,n\rangle\}\cup (C_\alpha\setminus F\times\{n-1\}):\alpha \text{ a limit }, n>0\}$ to be a weak neighborhood base. Of course, defining the space in this way will not be first countable or locally compact. But to see how to make such a space first countable and locally compact, see  \cite{VD}.

The $\clubsuit^\ast$ sequence $\{C_\alpha : \alpha\in\text{Lim}(\omega_1)\}$ and the convergence property of $C_\alpha\times\{n-1\}$ to $\langle\alpha,n\rangle$ are used to prove Lemma \ref{unctble} and Lemma \ref{cocount}, as well as to demonstrate that the space is not semi-proximal.

The usual de Caux space satisfies the version of the following with ``uncountable'' in place of ``club''.

For $A\subset X$, define $A(n) = \{\alpha: \langle\alpha,n\rangle\in {A}\}$, for all $n<\omega$.
\begin{lemma}\label{unctble} If $A\subset (\omega_1\times\{k\})$ is uncountable, then $\overline{A}(n)$ contains a club for all $n>k$\end{lemma} 
\begin{proof} We will prove it by induction on $n> k$, so it is sufficient to prove that if $\overline{A}$ is uncountable in $\omega_1\times\{n\}$, then $\overline{A}({n+1})$ is a club. Let $n> k$ and let $\overline{A}$ is uncountable in $\omega_1\times\{n\}$, then $|C_\alpha\times \{n\}\cap \overline{A}| = \omega$ for a club many $\alpha$. Hence, $\langle\alpha,n+1\rangle\in \overline{A}$, for a club many $\alpha$, so that $\overline{A}({n+1})$ is a club. \end{proof}

Therefore, if $A$ and $B$ are uncountable subsets of $X$, then there exist $i,j<\omega$ such that $A\cap (\omega_1\times \{i\})$ and $B\cap (\omega_1\times \{j\})$ are uncountable. By lemma \ref{unctble},  $\overline{A}(n)$ and $\overline{B}(n)$ are clubs for all $n>\max\{i,j\}$. Since $\overline{A}(n)\cap \overline{B}(n)$ is a club, then $\overline{A}\cap\overline{B}\cap(\omega_1\times\{n\})$ is uncountable for all $n> \max\{i,j\}$.  As a result, if $H$ and $K$ are disjoint closed subsets of $X$, then one of them must be countable.

To see that $X$ is normal, suppose that $H$ and $K$ are  disjoint closed subsets of $X$, then one of them, say $K$, is countable. Thus there is an $\alpha$ such that $K\subset (\alpha+1)\times \omega$. Since $K$ and $H^\prime = H\cap( (\alpha+1)\times \omega)$ are disjoint closed subset of a countable regular clopen subspace $(\alpha+1)\times \omega$, there exists two disjoint open set $U$ and $V^\prime$ in $(\alpha+1)\times \omega$ such that $K\subset U$ and $H^\prime\subset V^\prime$. Now let $V = V^\prime \cup ((\omega_1\setminus \alpha+1)\times \omega)$, then $U$ and $V$ are disjoint open set containing $K$ and $H$ respectively.

 $X$ is not countably paracompact since if $D_n = \omega_1\times (\omega\setminus n)$ and $U_n$ is open set containing $D_n$, for all $n<\omega$. Then $X\setminus U_n$ and $D_n$ are disjoint closed subsets of $X$. Thus $X\setminus U_n$ is countable and thus $\bigcap\{U_n : n<\omega\}\neq \emptyset$.

The usual de Caux space satisfies the version of the following lemma with the stronger assumption that $S$ is club. 

\begin{lemma}\label{cocount0}
    For $n\in\omega$, if $S$ is a stationary subset of $\omega_1$ and $U$ be any open neighborhood of $S\times\{n\}$, then $U\cap(\omega_1\times\{n-1\})$ is co-countable.
\end{lemma}
\begin{proof} Fix $n\in\omega$. Let $S$ be a stationary subset of $\omega_1$ and $U$ be any open neighborhood of $S\times\{n\}$. Suppose $A = (\omega_1\times\{n-1\})\setminus U$ is uncountable. Then there exists a club $C$ such that $A(n-1)\cap C_\alpha$ is infinite for all $\alpha\in C$. But for every $\alpha\in C\cap S$, $\langle\alpha, n\rangle\in U$ and hence $C_\alpha\times\{n-1\}\subset^*U$, contradiction! Thus, $U\cap(\omega_1\times\{n-1\})$ is co-countable.\end{proof}

\begin{lemma}\label{cocount} If $S$ is a stationary subset of $\omega_1$, then for every $n\in\omega$ and every $U$ open neighborhood of $S\times\{n\}$, $U\cap(\omega_1\times\{k\})$ is co-countable, for all $k<n$. 
\end{lemma}

\begin{proof}Let $S$ be a stationary subset of $\omega_1$.  We will prove it by induction on $n$.  If $n =1$, let $U$ be any open neighborhood of $S\times\{1\}$.  Then $U\cap (\omega_1\times\{0\})$ is co-countable, by Lemma \ref{cocount0}. Assume it is true for $n-1$, and let $U$ be any open neighborhood of $S\times\{n\}$. Then, by Lemma \ref{cocount0}, $U\cap (\omega_1\times\{n-1\})$ is co-countable. Thus,  there exists $\eta\in\omega_1$ such that $\omega_1\setminus\eta\times\{n-1\} \subset U$. Since $\omega_1\setminus\eta$ is stationary then,  by inductive hypothesis, $U\cap(\omega_1\times\{k\})$ is co-countable, for all $k<n-1$ and this concludes the proof. 
\end{proof}

\begin{lemma}\label{UniStDow} Fix $n\in \omega$. For every uniformity ${\mathfrak U}$ on $X$, $U\in{\mathfrak U}$, and $S$ a stationary subset of $\omega_1$, there exist $\beta\in S$ and a stationary set $S^\prime\subset S$ such that $S^\prime\times\{n+1\}\subset U[\langle\beta,n\rangle]$. \end{lemma}

\begin{proof}  Let $n\in\omega$ and $S$ be a stationary subset of $\omega_1$. Since $S\times\{n\}$ is uncountable, then $\overline{S\times\{n\}}(k)$ is a club for all $k>n$, by Lemma \ref{unctble}. Let ${\mathfrak U}$ be a uniformity on $X$ and let $U\in{\mathfrak U}$. Let $C = \overline{S\times\{n\}}(n+1)$. Then for each $\alpha\in C$, $A_\alpha = U[\langle\alpha, n+1\rangle]\cap S\times\{n\}\neq\emptyset$. 
Define $f: C \to S$ by $f(\alpha) = \min(A_\alpha(n))$.   Then by Fodor's lemma, there exist $\beta\in S$ and a stationary set $S^\prime\subset C\cap S$ such that $\beta= f(\alpha)$, for all $\alpha\in S^\prime$. Hence $S^\prime\times\{n+1\}\subset U[\langle\beta,n\rangle]$. \end{proof}

To see that $X$ is not semi-proximal, let ${\mathfrak U}$ be any uniformity on $X$. We will define a winning strategy $\sigma$ for Player II in the proximal game on $(X,{\mathfrak U})$. In inning $0$, Player I plays $U_0\in{\mathfrak U}$. Then, by Lemma \ref{UniStDow}, there exist $\beta_0\in\omega_1$ and a stationary set $S_0\subset \omega_1$ such that $S_0\times\{1\}\subset U[\langle \beta_0,0\rangle]$. In inning $n>0$, assume we have defined $U_m$, $S_m$, and $\beta_m$ for all $m<n$ such that:
\begin{enumerate}
\item $U_m$ is the response of Player I against $\sigma$.
\item $S_m\subset S_{m-1}$ for all $m>0$, and $S_m$ is stationary for all $m<n$. 
\item $\beta_m\in S_{m-1}\setminus(\beta_{m-1}+1)$, $\langle\beta_{m},m\rangle = \sigma(U_0,\cdots, U_m)$, and $S_m\times\{m+1\}\subset U[\langle\beta_{m},m\rangle]$.
\end{enumerate}

   In inning $n$, Player I plays $U_n\subset U_{n-1}$. Since $S_{n-1}$ is stationary, then there exist $\beta_n\in S_{n-1}\setminus(\beta_{n-1}+1)$ and a stationary set $S_n\subset S_{n-1}$ such that $S_n\times\{n+1\}\subset U_n[\langle\beta_n, n\rangle]$. And Player II chooses $\sigma(U_0,\cdots,U_n) = \langle \beta_n, n\rangle$. Note that, by inductive hypothesis, $\langle \beta_n, n\rangle\in U_{n-1}[\langle \beta_{n-1}, n-1\rangle]$ since $\beta_n\in S_{n-1}$. Therefore,  $(\langle \beta_n,n\rangle: n\in\omega)$ is not convergent in $X$ since the sets $W_n=\omega_1\times n$ form an increasing open cover of $X$ and $\langle\beta_k,k\rangle\not\in W_n$, for all $k\geq n$.  And since $S_{n}\times\{n+1\} \subset U_n[ \langle \beta_n, n\rangle]$ for all $n$, then by Lemma \ref{cocount},  $U_n[ \langle \beta_n, n\rangle]\cap \omega_1\times\{0\}$ is cocountable, for all $n$ and thus $\bigcap_{n\in\omega} U_n[ \langle \beta_n, n\rangle]\neq\emptyset$.

We remark that it seems rather delicate to determine if a Dowker space can be semi-proximal. Of course, any example that is not Fr\'echet is not, but we have no examples of semi-proximal Dowker spaces. Of course, a ZFC example is preferred but we do not even know of a consistent example: 

\begin{problem} Can there exist a semi-proximal Dowker space?
\end{problem}
\begin{problem}
    If $X$ is semi-proximal and normal space, is $X$ countably paracompact?
\end{problem}
\section{A semi-proximal non-normal $\Psi$-space}

Let $Z\subseteq 2^\omega$ and for each $z\in Z$,  Let $a_z = \{z\upharpoonright n : n\in\omega\}\subseteq 2^{<\omega}$. Then $A_Z = \{a_z:z\in Z\}$ is an almost disjoint family of branches in $2^{<\omega}$. Let $\Psi(A_Z) = 2^{<\omega}\cup A_Z$ with the points of $2^{<\omega}$ isolated and a local base at each $a_z$ is of the form $\{a_z\}\cup a_z\setminus F$ where $F$ is finite. For $s\in 2^{<\omega}$, let $[s]=\{t\in 2^{< \omega}: t\text{ end-extends }s\}$.

\begin{theorem}\label{psi} If $Z\subset 2^\omega$ contains no copy of the Cantor set, then $\Psi(A_Z)$ is semi-proximal.
\end{theorem} 
\begin{proof} Let $Z\subset 2^\omega$ contain no copy of the Cantor set and consider the uniformity ${\mathfrak U}$ induced on $\Psi(A_Z)$ as a subspace of the Stone-\v{C}ech compactification on $\Psi(A_Z)$. A base for this uniformity consists of entourages of the form $\bigcup\{U^2:U\in {\mathcal A}\}$ where ${\mathcal A}$ is a finite clopen partition of $\Psi((A_Z)$. In what follows, we will play the open cover version of the proximal game, where Player I will play finite clopen partitions of the space. 

\vspace{0.2cm}
   For $i = 0,1$, define \begin{align}\tag{*}
A_n^i  = \bigcup\{[s]: s\in 2^{n+1},s(n) = i\}\cup \{a_\alpha \in A_Z: z_\alpha(n) = i\}\label{*}
\end{align}
and for a finite subset $F\subset A_Z$ and $n\in\omega$, define the following:
\begin{enumerate} 
\item If $F = \{a\}$, ${\mathcal A}_{F,n} = \{\{a\}\cup a\setminus 2^{\leq n}\}$. 
\item If $|F|>1$,  fix $n_F\geq n$ so that $\{\{a\}\cup a\setminus (2^{\leq n_F}) : a\in F\}$ is pairwise disjoint, and define ${\mathcal A}_{F,n} = \{\{a\}\cup a\setminus 2^{\leq n_F} : a\in F\}$.

\item  $ {\mathcal A}_n = \{A_n^0, A_n^1\}\cup \{\{s\}: s\in 2^{\leq n} \}.$
\item 
${\mathcal A}_n\sim F  =   \{A_n^0\setminus \bigcup {\mathcal A}_{F,n}, A_n^1\setminus \bigcup {\mathcal A}_{F,n}\}\cup {\mathcal A}_{F,n}\cup \{\{s\}: s\in 2^{\leq n_F}\}.
$
\end{enumerate}
Note that for each $n$, ${\mathcal A}_n$ and ${\mathcal A}_n\sim F$ are both clopen partitions of the $\Psi$-space and hence Player I is free to play partitions of these forms in a play of the proximal game with respect to the uniformity ${\mathfrak U}$.

\vspace{0.2cm}
   Let $\sigma$ be a strategy for Player II in the proximal game on $(\Psi(A_Z), {\mathfrak U})$. 
Every play by Player I against $\sigma$ will be a clopen partition of the form ${\mathcal A}_n\sim F$. The {\em plain strategy} for Player I is to play ${\mathcal A}_n$ at every inning.  And a finite modification of the plain strategy for Player I is to play partitions of the form ${\mathcal A}_n\sim F$. 

\vspace{0.2cm}
   We will now define, for each $g\in 2^\omega$, plays of the game $P_g$ where Player II uses the strategy $\sigma$, so that if Player II wins each of these plays of the game, then a copy of the Cantor set would be embedded in $Z$. We first define a play of the game for the constant $0$ function:

\vspace{0.2cm}
   Let $F_{\langle \overline{0}\rangle\upharpoonright n+1} = \emptyset$, for all $n\geq 0$. In inning 0,  Player I chooses ${\mathcal A}_0\sim F_{\langle 0\rangle}$, and  Player II chooses $\sigma({\mathcal A}_0\sim F_{\langle 0\rangle}) = x_{{\langle 0\rangle}}$, which gives an initial play of the game denoted by $P_{\langle0\rangle}$.
Extend it to a full play of the game where Player I uses the unmodified plain strategy and Player II uses $\sigma$ to obtain 
 $$
 P_{\langle\overline{ 0}\rangle} = P_{\langle0\rangle} \smallfrown({ \mathcal A}_1\sim F_{\langle 00\rangle}, x_{\langle 00\rangle},\cdots,{ \mathcal A}_n\sim F_{\langle \overline{0}\rangle\upharpoonright n+1}, x_{\langle \overline{0}\rangle\upharpoonright n+1},\cdots).
 $$
 If either
\begin{enumerate}
    \item  there exists $n\geq 0$ such that $x_{\langle \overline{0}\rangle\upharpoonright n+1} \in 2^n$, or
    \item $\bigcap_{n\in\omega} \,\,\text{St}( x_{\langle \overline{0}\rangle\upharpoonright n+1}, {{\mathcal A}_n}) =\emptyset$
\end{enumerate} 
then $\sigma$ is defeated. Indeed, if (1) occurs, then Player II is forced to play $x_{\langle \overline{0}\rangle\upharpoonright n+1}$ in each subsequent inning and so picks a convergent sequence. So we assume that (1) and (2) fail. Then for each $n$ there is $i_{\langle \overline{0}\rangle\upharpoonright n+1}$ such that $x_{\langle \overline{0}\rangle\upharpoonright n+1} \in A^{i_{\langle \overline{0}\rangle\upharpoonright n+1}}_n$. Since the plain strategy was employed, there is $\alpha_{\langle{0}\rangle}$ such that $$\bigcap_{n\in\omega} \,\,\text{St}( x_{\langle \overline{0}\rangle\upharpoonright n+1}, {{\mathcal A}_n}) = \bigcap_{n\in\omega} A_n^{i_{\langle \overline{0}\rangle\upharpoonright n+1}} =
\{a_{\alpha_{\langle{0}\rangle}}\}.$$

\vspace{0.2cm}
   Now we use $a_{\alpha_{\langle{0}\rangle}}$ to define another play of the game corresponding to the branch $(1,0,0, \cdots)$ in $2^\omega$. Let  $F_{\langle 1\rangle} = F_{\langle 1\overline{0}\rangle\upharpoonright n+1} = \{a_{\alpha_{\langle 0\rangle}}\} $, for all $n> 0$.  
In inning 0,  Player I uses the plain strategy modified by $F_{\langle 1\rangle}$  and chooses $ { \mathcal A}_0\sim F_{\langle 1\rangle}$. While Player II chooses $x_{\langle 1\rangle}$. It gives an initial play of the game
$$
P_ {\langle 1\rangle} = ({{\mathcal A}_0\sim F_{\langle 1\rangle}},x_{\langle 1\rangle}).
$$
We extend it to a full play of the game with Player I using the plain strategy only modified by $F_ {\langle 1\rangle}$ and Player II using $\sigma$:
$$
P_ {\langle 1\overline{ 0}\rangle} =  P_{\langle1\rangle}\smallfrown({{ \mathcal A}_1\sim F_{\langle 10\rangle}}, x_{\langle 10\rangle},\cdots,{ \mathcal A}_n\sim F_{\langle 1\overline{0}\rangle\upharpoonright n+1}, x_{\langle 1\overline{0}\rangle\upharpoonright n+1},\cdots).
$$ 
 If either
\begin{enumerate}
    \item  there exists $n$ such that $x_{\langle 1\overline{0}\rangle\upharpoonright n+1}\in 2^{\leq n}$, or
    \item There exists $n$ such that $(x_{\langle 1\overline{0}\rangle\upharpoonright n+1})\in (a_{\alpha_{\langle0\rangle}}\setminus 2^{\leq n})\cup\{a_{\alpha_{\langle0\rangle}}\}$, or
    \item $\bigcap_{n\in\omega} \,\,\text{St}( x_{\langle 1\overline{0}\rangle\upharpoonright n+1}, {{\mathcal A}_n\sim F_{\langle1\overline{0}\rangle\upharpoonright n+1}}) =\emptyset$
\end{enumerate} 
Then, $\sigma$ is defeated.  Indeed, if (1) or (3) occurs, this is clear and if (2) occurs, then Player II is forced to play the rest of the game inside the convergent sequence $a_{\alpha_{\langle0\rangle}}\cup\{a_{\alpha_{\langle0\rangle}}\}$ and so Player I can play to defeat $\sigma$ (to produce either an eventually constant sequence or one that converges to $a_{\alpha_{\langle0\rangle}}$). 
Otherwise, there exists $\alpha_{\langle 1\rangle}\in\omega_1$ such that 
$$
\bigcap_{n\in\omega} \,\,\text{St}( x_{\langle 1\overline{0}\rangle\upharpoonright n+1}, {{\mathcal A}_n\sim F_{\langle1\overline{0}\rangle\upharpoonright n+1}}) = \bigcap_{n\in\omega} (A_n^{i_{\langle 1\overline{0}\rangle\upharpoonright n+1}}\setminus \bigcup {\mathcal A}_{F_{\langle 1\rangle},n}) = 
\{a_{\alpha_{\langle{1}\rangle}}\}.$$  
Note that $\alpha_{\langle 1\rangle}\neq \alpha_{\langle 0\rangle}$, so there exists a minimum $m_\emptyset\in\omega$ such that, in inning $m_\emptyset$,  $i_{\langle \overline{0}\rangle\upharpoonright m_\emptyset+1}\neq i_{\langle 1\overline{0}\rangle\rangle\upharpoonright m_\emptyset+1}$ which means that $z_{\alpha_{\langle 0\rangle}}( m_\emptyset) \neq z_{\alpha_{\langle 1\rangle}}( m_\emptyset)$ and $z_{\alpha_{\langle 0\rangle}}\upharpoonright m_\emptyset = z_{\alpha_{\langle 1\rangle}}\upharpoonright m_\emptyset$.

\vspace{0.2cm}
   Let $n>1$ and assume we have defined $P_s$, $\alpha_s$ and $F_{s}$, for every $s\in 2^{\leq n}$,  $m_t\in\omega$ for every $t\in 2^{< n}$, and we have also defined $F_{s\smallfrown\overline{0}\upharpoonright k}$ and  $i_{s\smallfrown\langle \overline{0}\rangle\upharpoonright k}$ for all $k>|s|$ such that:
\begin{enumerate} 
\item $P_s =({{\mathcal A}_0\sim F_{s \upharpoonright 1}}, x_{s \upharpoonright 1},\cdots, {{\mathcal A}_{n-1}\sim F_{s\upharpoonright n}}, x_{s\upharpoonright n })$ which is an initial play of the game at inning $n-1$. If $s$ extends $t$, then $P_s$ extends $P_t$.

$P_{s\smallfrown\langle \overline{0}\rangle} = P_s\smallfrown ({{\mathcal A}_n\sim F_s}, x_{s\smallfrown\langle \overline{0}\rangle\upharpoonright n+1},{{\mathcal A}_{n+1}\sim F_s}, x_{s\smallfrown\langle \overline{0}\rangle\upharpoonright n+2},\cdots)$ which is a branch of play of the game corresponds to $s\smallfrown\langle\overline{0}\rangle$. 

\item For all $s\in 2^{\leq n}$, $F_s = \{a_{\alpha_t} : s\in[t]\}$. And for $k>|s|$,  $F_{s\smallfrown\overline{0}\upharpoonright k} = F_s$.

\item $i_{s\smallfrown\langle \overline{0}\rangle\upharpoonright k+1}\in 2$ such that $x_{s\smallfrown \langle\overline{0}\rangle\upharpoonright k+1}\in  (A_k^{i_{s\smallfrown\langle \overline{0}\rangle\upharpoonright k+1}}\setminus \bigcup {\mathcal A}_{F_{s\smallfrown\overline{0}\upharpoonright k+1}, k})$, for all $k$.

\item $\{a_{\alpha_s}\} = \bigcap_{k\in\omega} (A_k^{i_{s\smallfrown\langle \overline{0}\rangle\upharpoonright k+1}}\setminus \bigcup {\mathcal A}_{F_{s\smallfrown\overline{0}\upharpoonright k+1}, k}$).

\item  $m_t$ satisfies the following:
\begin{enumerate}
    \item $z_{\alpha_{t\smallfrown 0}} ( m_t) \neq z_{\alpha_{t\smallfrown 1}} ( m_t )$, and  $z_{\alpha_{t\smallfrown 0}} \upharpoonright m_t = z_{\alpha_{t\smallfrown 1}} \upharpoonright  m_t $. 
    \item  $z_{\alpha_r}\upharpoonright m_r = z_{\alpha_t}\upharpoonright m_r$ if $t\in[r]\cap2^{<n}$. 
    \item $z_{\alpha_r}\upharpoonright m_{\Delta(r,t)} = z_{\alpha_t}\upharpoonright m_{\Delta(r,t)} \,\,\textnormal{and}\,\, z_{\alpha_r}(m_{\Delta(r,t)}) \neq z_{\alpha_t}(m_{\Delta(r,t)})$ if $r$ and $t$ are not comparable in $2^{<n}$, where $\Delta(r,t)$ is the maximum in $2^{<n}$ which both $r$ and $t$ extend.
\end{enumerate}
\end{enumerate}
 Now, for every $s\in 2^{n}$, define $F_{s\smallfrown 0} =  F_s$, and $F_{s\smallfrown0\smallfrown\langle\overline 0 \rangle\upharpoonright k+1} = F_s$, for all $k> n$.  
 In inning $n$, the initial play of the game, where Player I uses the plain strategy modified by $F_{s}$ and Player II uses $\sigma$, is $P_{s\smallfrown 0} = P_s\smallfrown ({{\mathcal A}_n\sim F_{s}}, x_{s\smallfrown 0})$, which we extend to a full play of the game with the plain strategy only modified by $F_{s}$,
 $$
 P_{s\smallfrown\langle0 \overline{0}\rangle} = P_{s\smallfrown 0}\smallfrown ({{\mathcal A}_{n+1}\sim F_s}, x_{s\smallfrown\langle \overline{0}\rangle\upharpoonright n+2}\cdots).
 $$
  Therefore, the play of the game, corresponds to $s\smallfrown 0$, is equal to the one that corresponds to $s$, $P_{s\smallfrown\langle \overline{0}\rangle}$. So, let $\alpha_{s\smallfrown 0} = \alpha_s$.

  \vspace{0.2cm}
    Define $F_{s\smallfrown 1} =  F_s\cup\{a_{\alpha_s}\}$ and
$F_{s\smallfrown\langle1\overline{0}\rangle\upharpoonright k+1} = F_{s\smallfrown 1}$, for all $k>n$. In inning $n$, the initial play of the game, where Player I uses the plain strategy modified by $F_{s\smallfrown 1\upharpoonright k+1}$ at inning $k\leq n$ and Player II uses $\sigma$, is $P_{s\smallfrown 1} = P_s\smallfrown ({{\mathcal A}_n\sim F_{s\smallfrown 1}}, x_{s\smallfrown 1} )$, which we extend to a full play of the game with the plain strategy modified by $F_{s\smallfrown 1}$ at inning $k>n$, 
 $$
 P_{s\smallfrown\langle1 \overline{0}\rangle} = P_{s\smallfrown 1}\smallfrown ({{\mathcal A}_{n+1}\sim F_{s\smallfrown 1}}, x_{s\smallfrown\langle1 \overline{0}\rangle\upharpoonright n+2},\cdots).
 $$
As in the base case of the construction, if either
\begin{enumerate}
    \item there exists $k$ such that $x_{s\smallfrown\langle 1\overline{0}\rangle\upharpoonright k+1}\in 2^{ k}$ 
    \item there exists $a\in F_{s\smallfrown 1}$ such that $(x_{s\smallfrown\langle 1\overline{0}\rangle\upharpoonright k+1})\in (a\setminus{2^{\leq n_{F_{s\smallfrown 1}}}})\cup\{a\}$, or
    \item $\bigcap_{k\in\omega} \,\,\text{St}( x_{s\smallfrown\langle 1\overline{0}\rangle\upharpoonright k+1}, {{\mathcal A}_k\sim F_{s\smallfrown\langle 1\overline{0}\rangle\upharpoonright k+1}}) =\emptyset$
\end{enumerate} 
Then, there is a play of the game where the strategy $\sigma$ is defeated. Otherwise, for each $s\in 2^n$ there exists $\alpha_{s\smallfrown 1}\in\omega_1$ such that 
$$
\bigcap_{k\in\omega} \,\,\text{St}( x_{s\smallfrown\langle1 \overline{0}\rangle\upharpoonright k+1}, {\mathcal A}_k\sim F_{s\smallfrown\langle 1\overline{0}\rangle\upharpoonright k+1}) = \bigcap_{k\in\omega}( A_n^{i_{s\smallfrown\langle 1\overline{0}\rangle\upharpoonright k+1}}\setminus \bigcup {\mathcal A}_{F_{s\smallfrown\langle 1\overline{0}\rangle\upharpoonright k+1}, k}) = \{a_{\alpha_{s\smallfrown 1}}\}. 
$$ 
Since $a_{\alpha_s}\in F_{s\smallfrown 1}$, then  $\alpha_{s\smallfrown 0} = \alpha_s  \neq \alpha_{s\smallfrown 1}$, so there exists a minimum $m_{s}$ such that in inning $m_s$, $i_{s\smallfrown\langle  \overline{0}\rangle\upharpoonright m_{s}+1}\neq i_{s\smallfrown \langle 1\overline{0}\rangle\upharpoonright m_{s}+1}$. Which means that $z_{\alpha_s}( m_s) \neq z_{\alpha_{s\smallfrown 1}}( m_s)$ and $z_{\alpha_{s}}\upharpoonright m_s = z_{\alpha_{s\smallfrown 1}}\upharpoonright m_s$.
 Then, the elements of $\{\alpha_s: s\in 2^n\}$ are distinct.  To see that, let $r\neq s$ in $2^n$. Consider $\Delta(s,r) = t\in 2^{k}$, for some $k<n$ which means that $s(k)\neq r(k)$. Assume without loss of generality that $s(k) = 0$ and $r(k) = 1$. Then we have, $$z_{\alpha_s}(m_t) = z_{\alpha_{t\smallfrown 0}}(m_t)  \neq z_{\alpha_{t\smallfrown 1}}(m_t)  = z_{\alpha_r}(m_t).$$

\vspace{0.2cm}
   If there is an $f\in 2^\omega$ such that 
$$\bigcap_{k\in\omega} (A_k^{i_{f\upharpoonright k+1}}\setminus \bigcup {\mathcal A}_{F_{f\upharpoonright k+1}, k}) = \emptyset,$$
then the play of the game, corresponding to that $f$, is the one when Player II uses $\sigma$, but it is winning for Player I. Otherwise for each $f\in 2^\omega$ there is a unique $\alpha_f$ such that $$\bigcap_{k\in\omega} (A_k^{i_{f\upharpoonright k+1}}\setminus \bigcup {\mathcal A}_{F_{f\upharpoonright k+1}, k}) = \{a_{\alpha_f}\}.$$ Hence, the mapping $G:2^\omega\to Z$ is one-to-one, where  $G(f) = z_{\alpha_f}$, for each $f\in 2^\omega$. It is also continuous. Indeed, let $V$ be an open set in $Z$, then there exists $s\in 2^n$ such that $V = [s]\cap Z$. Then there exists $t_s \in 2^n$ such that $s(k) = i_{t_s\upharpoonright k+1}$. Hence $G^{-1}(V) = [t_s]$ which is open in $2^\omega$.  Then $Z = \{z_\alpha : \alpha\in\omega_1\}$ contains a Cantor set. Contradiction.
\end{proof}

The following is a counterexample to Nyikos's question of whether normal implies semi-proximal. There will be a second example in section 4 and it  is a subspace of $(\omega+1)\times \omega_1$.  

\begin{example} There is a semi-proximal not normal $\Psi$-space. \end{example}

\begin{proof} Let $Z\subseteq 2^\omega$ be a Bernstein set. Consider the $\Psi$-space, $\Psi(A_Z)$. Since $Z$ has size continuum then $\Psi(A_Z)$ is not normal (by Jones Lemma since it is separable with a closed discrete set of size continuum). And since $Z$ contains no copy of the Cantor set, then $\Psi(A_Z)$ is semi-proximal, by Theorem \ref{psi}.
\end{proof}
\vspace{0.7cm}
 We can derive a corollary from our findings that relates to the Galvin game, as described in \cite{MS} (also called the weak Ulam game in \cite{Ph}). In this game, Player I starts by playing a finite partition ${\mathcal A}_0$ of $X$, and then Player II chooses an element $A_0$ from ${\mathcal A}_0$ with $|A_0|>1$. In each subsequent inning, $n > 0$, Player I plays a finite partition ${\mathcal A}_n$ of $X$ refining ${\mathcal A}_{n-1}$ (or equivalently a finite partition of the set chosen by Player II in the previous inning), and Player II responds by picking an element $A_n\in{\mathcal A}_n$ with $A_n\subseteq A_{n-1}$ with $|A_n|>1$. Player I wins and the game immediately ends if Player II ever picks a singleton. Player II wins the play of the game if the intersection of the sets $\{A_n : n < \omega\}$ is non-empty. The following corollary is due to Galvin (\cite{MS}) but follows from our Theorem \ref{psi}.
\begin{corollary} Player II does not have a winning strategy in the Galvin game on $\kappa$ for $\kappa\leq 2^{\aleph_0}$. 
\end{corollary} 
\begin{proof} Identify $\kappa$ with a subset $Z$ of the Cantor set that contains no copies of the Cantor set. For example, in the case $\kappa=\mathfrak{c}$, let $Z$ be a Bernstein set. Let $\sigma$ be a strategy for Player II in the Galvin game on $\kappa$. If Player II ever chooses a finite set, then Player I can defeat it by partitioning the set into singletons. So, assume Player II always chooses an infinite set using $\sigma$.  

\vspace{0.2cm}
   Consider the uniformity ${\mathfrak U}$ on $\Psi(A_Z)$ inherited from the Stone-\v{C}ech compactification of $\Psi(A_Z)$. Define a strategy $\sigma^\prime$ for Player II in the proximal game on $(\Psi(A_Z), {\mathfrak U})$ as follows: For a finite sequence of finite clopen partitions ${\mathcal A}_0, \cdots, {\mathcal A}_n$ of $\Psi(A_Z)$ such that ${\mathcal A}_n$ refines ${\mathcal A}_{n-1}$, let $\sigma^\prime({{\mathcal A}_0},\cdots, {{\mathcal A}_n}) =  a_{\alpha_n}$ such that $\alpha_n\in \sigma({\mathcal A}^\prime_0,\cdots, {\mathcal A}^\prime_n)$ and $a_{\alpha_n}\notin\{a_{\alpha_k}:k<n\}$ where $${\mathcal A}^\prime_k = \{\{\alpha\in \kappa: a_\alpha\in U\cap A_Z\}: U\in{\mathcal A}_k\},$$ for each $k\leq n$. 

\vspace{0.2cm}
   Since $(\Psi(A_Z),{\mathfrak U})$ is semi-proximal by Theorem \ref{psi} and since $(a_{\alpha_n})_{n \in \omega}$ does not converge in $\Psi(A_Z)$ as it is a non eventually constant sequence in the discrete closed set $A_Z$, there exists a sequence $({ \mathcal A}_n: n\in\omega)$ of finite clopen partition of $\Psi(A_Z)$ such that $\bigcap_{n\in\omega}\,\,\text{St}(a_{\alpha_n},{\mathcal A}_n) =\emptyset$. Hence, $\bigcap_{n\in\omega}\sigma({\mathcal A}^\prime_0,\cdots, {\mathcal A}^\prime_n)=\emptyset$. Thus, $\sigma$ is not winning. 
\end{proof}


\section{Semi-proximal and normal subspaces of finite powers of $\omega_1$}

 Fleissner proved the following theorem in \cite{Fl}: 
\begin{theorem}\label{Fl} The following are equivalent for a subspace $X$ of a finite products of ordinals: 
\begin{enumerate}
    \item $X$ is normal.
    \item $X$ is normal and strongly zero-dimensional.
    \item $X$ is collectionwise normal.
    \item $X$ is shrinking. 
\end{enumerate}
\end{theorem} 
It is natural to ask what can be said in general about semi-proximality of subspaces of products of ordinals. We can restrict our study to powers of $\omega_1$ since subspaces of $\kappa$ for  $\kappa>\omega_1$ are not in general Fr\'echet, so not semi-proximal (even if normal). So we consider the relationship between normality and semi-proximality in subspaces of finite products of $\omega_1$.

  For a finite ordinal $k$ and for $X\subseteq \omega_1^k$, fix the following notations:
\begin{enumerate}
    \item $\Delta(X)=\{\alpha:\text{the $k$-tuple }(\alpha,....,\alpha)\in X\}$.
    \item For $\alpha\in\omega_1$,  $X_\alpha = \{ \langle \delta_0,\cdots,\delta_{k-1}\rangle\in\omega_1^{k-1}: \langle \alpha,\delta_0,\cdots,\delta_{k-1}\rangle \in X\}.$ 
    \item For $j<k$ and $\alpha\in\omega_1$, let $X_\alpha^{j} = \{\delta\in\omega_1:\exists x\in X_\alpha (\delta = \delta_j)\}.$
    \item For $A,B\subset \omega_1$, let 
    $X^k_A  =X \cap(A\times\omega_1^k)$ and $X_A^{B^k} = X \cap( A \times B^k)$ .
\end{enumerate}
\begin{theorem} \label{NormSemi} Every normal subspace of finite powers of $\omega_1$ is semi-proximal. \end{theorem} 

\begin{proof}
    Let $k$ be a finite ordinal and let $X$ be a normal subspace of $\omega_1^k$. If $k=1$, then it is true since $\omega_1$ is hereditarily semi-proximal by \cite{Pa}. Now, assume the theorem is true for all $j<k$.
\begin{claim} \label{SemiNorm}
      For $\alpha\in\omega_1$,  $X^{k-1}_{[0,\alpha]}$ is semi-proximal.
\end{claim}

\begin{proof}
    We will prove it by induction on $\alpha$. For $\alpha = 1$,  $X$ is a normal subspace of a homeomorphic to $\omega_1^{k-1}$, $\{0\}\times \omega_1^{k-1}$. Hence, $X$ is semi-proximal by the induction hypothesis on $k$.  Now, assume it is true for all $\gamma <\alpha$.  If $\alpha$ is a successor, then $X$ is semi-proximal since it is a disjoint union of clopen semi-proximal subspaces. If $\alpha$ is a limit ordinal, we have two cases:

    \vspace{0.2cm}
   Case 1: There exists $j<k$ such that $X_\alpha^j$ is not stationary. In this case a club $C$ exists such that $C^{k-1}\cap X_\alpha= \emptyset$. Note that $X_{[0,\alpha)}^ {C^{k-1}}$ and $ X^{k-1}_{\{\alpha\}}$ are disjoint closed subsets of $X$. Given that $X$ is normal and strongly zero dimensional by Theorem \ref{Fl}, there exists a clopen set $U$ in $X$ such that $X_{[0,\alpha)}^ {C^{k-1}}\subseteq U$ and $X_{\{\alpha\}}^{k-1}\subseteq X\setminus U$.   Note that $U$ is normal since it is closed in $X$. We fix a strictly increasing cofinal sequence $\{\alpha_n: n \in \omega\}$ in $\alpha$. Since $ U = U^{k-1}_{[0,\alpha)}$ is represented as the direct sum $\bigoplus_{n\in\omega}U^{k-1}_{(\alpha_{n-1},\alpha_n]}$ of clopen subspace of $U$, where $\alpha_{-1} = -1$, then $U$ is semi-proximal by the inductive hypothesis on $\alpha$. Since each $U^{k-1}_{(\alpha_{n-1},\alpha_n]}$ is semi-proximal by the inductive hypothesis on $\alpha$.  Then, $U$ is semi-proximal. And $X\setminus U$ is semi-proximal since it is a subset of a metrizable space $(\alpha+1)\times(\omega_1^{k-1}\setminus C^{k-1})$. Therefore, $X$ is semi-proximal since it is a disjoint union of clopen semi-proximal subspaces.

\vspace{0.2cm}
   Case 2:  If $X_\alpha^j$ is stationary, for all $j<k$. Let ${\mathfrak U}$ be the uniformity induced on $X$ as a subspace of the compact space $(\alpha+1)\times (\omega_1+1)^{k-1}$. Therefore,  we will play the open cover version of the proximal game, where Player I will play finite clopen partitions of the space.  Let $\sigma$ be a strategy for Player II in the proximal game on $(X,{\mathfrak U})$. Let ${\mathcal M}$ be a countable elementary sub-model of a large enough portion of the universe containing $\sigma$ and $\bigcap_{j<k}X^j_\alpha$ such that $\lambda={\mathcal M}\cap \omega_1 \in \bigcap_{j<k}X_\alpha^j$. This can be done since $X_\alpha^j$ is stationary, for all $j<k$.   Let $\eta<\alpha$ and $\gamma\in\omega_1$, consider the following finite partition of $X$:
 $$
 {\mathcal U}_{\eta,\gamma} = \{X_{[0,\eta]}^{k-1}, X_{(\eta,\alpha]}^ {[0,\gamma]^{k-1}}, X_{(\eta,\alpha]}^ {(\gamma,\omega_1)^{k-1}}\}  
 $$
 Let $\{\lambda_n:n\in\omega\}$ be a strictly increasing sequence that converges to $\lambda$. Let  $x_n = \langle \beta_n, \delta^0_n,\cdots,\delta_n^{k-1}\rangle\in X$ be the point chosen by Player II during the nth inning, where $n\in\omega$, following the strategy $\sigma$. We define a winning counter-strategy $\tau_\sigma$ for Player I. In inning $n = 0,1$, let $\eta_n = \alpha_n$ and $\gamma_n = \lambda_n$. Define $\tau_\sigma(\emptyset) = { \mathcal U}_0$ and $\tau_\sigma(x_0) = { \mathcal U}_1$, where ${ \mathcal U}_0 = {\mathcal U}_{\eta_0,\gamma_0}$ and ${ \mathcal U}_1$ is the common refinement of  ${\mathcal U}_{\eta_0,\gamma_0}$ and ${\mathcal U}_{\eta_1,\gamma_1}$. Let $n>0$ and assume that, in inning $n+1$, we have chosen $\eta_{n+1}$ such that $\eta_n<\eta_{n+1}<\alpha$ and $\gamma_{n+1}$ such that $\gamma_n<\gamma_{n+1}<\omega_1$. And assume we have defined $\tau_\sigma(x_0,\cdots, x_n) = { \mathcal U}_{n+1}$, where ${ \mathcal U}_{n+1}$ is the common refinement of ${\mathcal U}_{\eta_0,\gamma_0},\cdots, {\mathcal U}_{\eta_{n+1},\gamma_{n+1}}$. Note that, $x_{n+1}\in U$, for some $U\in{ \mathcal U}_n$. If $U$ is not a subset of $X_{(\eta_n,\alpha]}^ {(\gamma_n,\omega_1)^{k-1}}$, then the remainder of the game would be conducted within a semi-proximal space. Consequently, there exists a play of the game where Player I could find a strategy to defeat $\sigma$. 
  Otherwise, $U_n\subset X_{(\eta_n,\alpha]}^ {(\gamma_n,\omega_1)^{k-1}}$, then let $\eta_{n+2}<\alpha$ such that $\eta_{n+2}>\max\{\eta_{n+1},\beta_{n+1}, \alpha_{n+2}\}$ and $\gamma_{n+2}\in{\mathcal M}$ such that $\gamma_{n+2}> \max\{ \gamma_{n+1}, \lambda_{n+2}\}\cup\{\delta_{n+1}^j:j<k\}$, and define $\tau_{\sigma}(x_0,\cdots,x_{n+1}) = { \mathcal U}_{n+2}$ such that ${\mathcal U}_{n+2}$ is the common clopen refinement of ${\mathcal U}_{\eta_0,\gamma_0},\cdots, {\mathcal U}_{\eta_{n+2},\gamma_{n+2}}$.
   If this is the case for all $n\in\omega$, then we get $\{ \langle\beta_n,\delta^0_n,\cdots,\delta^{k-1}_n\rangle\colon n\in\omega\}$ converges to $\langle \alpha,\lambda,\cdots,\lambda\rangle\in X$ and hence $\tau_\sigma$ is winning. Hence, $X$ is semi-proximal.
\end{proof}

  Now, consider two cases for $X\subset \omega_1^k$:

Case 1: if $\Delta(X)$ is not stationary, then by (Lemma 3.5, \cite{Fl}), there exists a club $C$ on $\omega_1$ such that $X\cap C^k = \emptyset$. Let $V_i = \{ \langle \beta_0,\cdots,\beta_k\rangle\in\omega_1^k: \beta_i\notin C\}$. Note that each $V_i\cap X$ can be written as a direct sum of $\bigoplus_{\alpha\in C}X_{(\alpha,\alpha^+]}^{k-1}$, where $\alpha^+ = \min C\setminus\alpha$. Therefore, each $V_i\cap X$ is semi-proximal by claim \ref{SemiNorm}. By Theorem \ref{Fl}, $X$ is normal and strongly zero-dimensional. As a result, the finite open cover $\{V_i\cap X : i\leq k\}$ of $X$ has a disjoint clopen refinement $\{W_i:i\leq k\}$. Each $W_i$ is semi-proximal, implying that $X$ is semi-proximal.

   Case 2:   $\Delta(X)$ is stationary. Let ${\mathfrak U}$ be the uniformity induced on $X$ as a subspace of $(\omega_1+1)^k$. Therefore,  we will play the open cover version of the proximal game, where Player I will play finite clopen partitions of the space.  Let $\sigma$ be a strategy for Player II in the proximal game on $(X,{\mathfrak U})$, and let ${\mathcal M}$ be a countable elementary submodel of a large enough portion of the universe containing $\sigma$ and $\Delta(X)$. Since $\Delta(X)$ is stationary, let $\lambda = {\mathcal M}\cap \omega_1\in\Delta(X)$. Let $\{\lambda_n:n\in\omega\}$ be a strictly increasing sequence that converges to $\lambda$. For every $\gamma\in\omega_1$, consider the following finite partition of $X$:
$$
{\mathcal U}_\gamma = \{\prod_{i\leq k} A_i\cap X : A_i\in\{(\gamma,\omega_1), [0,\gamma]\}\}
$$
 Let $x_n = \langle \delta_0^n,\cdots, \delta_k^n\rangle$ be the point chosen by Player II using $\sigma$. We will define a winning counter-strategy $\tau_\sigma$ for Player I.  In inning $n = 0,1$, let $\gamma_n = \lambda_n$. Define $\tau_\sigma(\emptyset) = { \mathcal U}_0$ and $\tau_\sigma(x_0) = { \mathcal U}_1$, where ${ \mathcal U}_0 = {\mathcal U}_{\gamma_0}$ and ${ \mathcal U}_1$ is the common refinement of  ${\mathcal U}_{\gamma_0}$ and ${\mathcal U}_{\gamma_1}$. Let $n>0$ and assume that, in inning $n+1$, we have chosen $\gamma_{n+1}$ such that $\gamma_n<\gamma_{n+1}<\omega_1$. And assume we have defined $\tau_\sigma(x_0,\cdots, x_n) = { \mathcal U}_{n+1}$, where ${ \mathcal U}_{n+1}$ is the common refinement of ${\mathcal U}_{\gamma_0},\cdots, {\mathcal U}_{\gamma_{n+1}}$. Note that, $x_{n+1}\in U$, for some $U\in{ \mathcal U}_n$. If $U$ is not a subset of $(\gamma_n,\omega_1)^k\cap X$, then the rest of the game would be played inside a semi-proximal space. Thus, there is a subsequent play of the game where Player I would find a way to defeat $\sigma$. 
  Otherwise, $U_n\subset (\gamma_n,\omega_1)^k\cap X$, then let  $\gamma_{n+2}\in{\mathcal M}$ such that $\gamma_{n+2}> \max\{ \gamma_{n+1}, \lambda_{n+2}\}\cup\{\delta_{n+1}^j:j<k\}$, and define $\tau_{\sigma}(x_0,\cdots,x_{n+1}) = { \mathcal U}_{n+2}$ such that ${\mathcal U}_{n+2}$ is the common clopen refinement of ${\mathcal U}_{\gamma_0},\cdots, {\mathcal U}_{\gamma_{n+2}}$.
   If this is the case for all $n\in\omega$, then we get $\{ \langle\delta^0_n,\cdots,\delta^{k}_n\rangle\colon n\in\omega\}$ converges to $\langle\lambda,\cdots,\lambda\rangle\in X$ and hence $\tau_\sigma$ is winning. Hence, $X$ is semi-proximal.
\end{proof}

The following is an example to show that the converse to Theorem \ref{SemiNorm} does not hold:

\begin{example}  There is a semi-proximal not normal subspace of $(\omega+1) \times \omega_1$. \end{example}
\begin{proof}
 Let $A = \{a_\alpha:  \alpha\in \text{Succ}\}$ be a family of subsets of $\omega$. Define $X_A\subseteq(\omega+1)\times\omega_1$ by $X_A = \big(\omega\times \text{Lim}\big)\cup(\bigcup_{\alpha\in\omega_1}(a_\alpha\cup\{\omega\}) \times \{\alpha+1\}).$ 

\vspace{0.2cm}
   This type of subspace was first introduced in \cite{Kem} by N. Kemoto where he proved that it is always not normal for any family $A$ as $H = \omega\times\text{Lim}$ and $K = \{\omega\}\times \text{Succ}$ cannot be separated. 

\vspace{0.2cm}
   Now, we prove that it is semi-proximal if $A$ is an almost disjoint family of the form ${A}_Z$ where $Z\subseteq 2^\omega$ contains no copy of the Cantor set.
 Let $Z = \{z_\alpha:\alpha\in\omega_1\}$ be a subset of $2^\omega$ and for each $\alpha\in \omega_1$,  let $a_\alpha = \{z_\alpha\upharpoonright n : n\in\omega\}\subseteq 2^{<\omega}$. Then $A_Z = \{a_\alpha:\alpha\in \omega_1\}$ is an almost disjoint family of branches in $2^{<\omega}$. Enumerate $2^{<\omega}$ as $\{s_n:n\in\omega\}$ and define $X_{A_Z} \subset (\omega+1)\times\omega_1$ by $$X_{A_Z} = \big(\omega\times \text{Lim}\big)\cup \big(\{\omega\}\times\text{Succ}\big)\cup\bigcup_{\alpha\in\omega_1}\{\langle n,\alpha+1\rangle:s_n\in a_\alpha\}.$$ 

 \vspace{0.2cm}
   To see that $X_{A_Z}$ is semi-proximal, consider the uniformity ${\mathfrak V}$ inherited from the Stone-\v{C}ech compactification of $X_{A_Z}$. So any clopen partition of $X_{A_Z}$ corresponds to an entourage in ${\mathfrak V}$ so we can consider the version of the game where Player I plays finite clopen partitions of the space. Let $\sigma$ be a strategy for Player II in the proximal game on $(X_{A_Z},{\mathfrak V})$.  Consider the uniformity ${\mathfrak U}$ on $\Psi(A_Z)$ inherited from the Stone-\v{C}ech compactification of $\Psi(A_Z)$. 
 
\vspace{0.2cm}
    Let $X_k = (\{k\}\times\omega_1)\cap X_{A_Z}$. For a subset $U$ of $\Psi(A_Z)$, define  $$U^\prime = \bigcup\{X_k: s_k\in U\}\cup \{\langle\omega,\alpha+1\rangle: a_\alpha\in U\}.$$ 
Note that if ${\mathcal U}$ is a finite clopen partition of $\Psi(A_Z)$, then ${\mathcal U}^\prime= \{U^\prime:  U\in{\mathcal U}\}$ is a finite clopen partition of $X_{A_Z}$ and so corresponds to an entourage in  ${\mathfrak V}$ so can be played by Player I in the proximal game.

\vspace{0.2cm}
   For $n\in\omega$ and $F\subset A_Z$, let $\alpha_{F} = \max\{\beta:a_\beta\in F\}$ and define 
$$
{\mathcal U}_{n,F} = \{X_{A_Z}\cap ((n+1)\times \omega_1), X_{A_Z}\cap((n,\omega]\times(\alpha_F+1)), X_{A_Z}\cap((n,\omega]\times (\alpha_F,\omega_1)) \}.
$$
 Now consider the clopen partition ${\mathcal A}_{n}\sim F$ and $n_F$ which are defined in the proof of theorem \ref{psi}, then 
  ${\mathcal A}'_{n,F} = \{U'\cap A : U\in{\mathcal A}_{n}\sim F, A\in{\mathcal U}_{n_F,F}\}$ is a finite clopen partition of $X_{A_Z}$. 

 \vspace{0.2cm}
   Define a strategy $\sigma^\prime$ for Player II in the proximal game on $(\Psi(A_Z), {\mathfrak U})$ as follows: For a finite sequence of clopen partitions ${\mathcal U}_0, \cdots, {\mathcal U}_n$ of $\Psi(A_Z)$ such that ${\mathcal U}_n$ refines ${\mathcal U}_{n-1}$, 
 $$\sigma^\prime({{\mathcal U}_0},\cdots, {{\mathcal U}_n}) = \begin{cases}
     a_{\alpha_n}  & \text{if }  \sigma({{\mathcal V}_0^\prime, \cdots, {\mathcal V}_n^\prime}) = \langle\omega,\alpha_n+1\rangle \\
  s_{k_n}  & \text{if }  \sigma({{\mathcal V}_0^\prime, \cdots, {\mathcal V}_n^\prime}) \in X_{k_n}
\end{cases}
$$
where ${\mathcal V}'_n $ is a finite clopen partition of $X_{A_Z}$ and the common refinement of ${ \mathcal U}_k'$, for all $k\leq n$. 
Since $(\Psi(A_Z),{\mathfrak U})$ is semi-proximal, then there is a play of the game where Player I defeats $\sigma^\prime$. 

\vspace{0.2cm}
   We need to recall some of the details from that proof. Recall that there is a sequence $(F_n:n\in\omega)$ of finite subsets of $A_Z$ such that the play of the game $$({\mathcal A}_0\sim F_0,x_0,\cdots,{\mathcal A}_n\sim F_n, x_n,\cdots)$$ is a winning play against $\sigma'$ and resulted in one of three outcomes:
\begin{enumerate}[(i)]
\item At some stage of the game, $({\mathcal A}_0\sim F_0,x_0,\cdots,{\mathcal A}_n\sim F_n,x_n)$ there is an $s\in 2^{\leq n_{F_n}}$ such that $\sigma'({\mathcal A}_0\sim F_0,x_0,\cdots,{\mathcal A}_n\sim F_n)=x_n=s$, or
\item At some stage of the game, $({\mathcal A}_0\sim F_0,x_0,\cdots,{\mathcal A}_n\sim F_n,x_n)$ there is an $a\in F_n$ such that $x_n\in \{a\}\cup a\setminus 2^{\leq n_{F_n}}$, or
\item The play of the game satisfied $\bigcap_{n\in\omega} \,\,\text{St}(x_n,{\mathcal A}_n\sim F_n) =\emptyset$.
\end{enumerate}
We consider the corresponding play of the proximal game on $X_{A_Z}$ given by the sequence of plays $\{{\mathcal V}'_n:n\in\omega\}$ by Player I against $\sigma$, where ${\mathcal V}'_n$ is the finite clopen partition of $X_{A_Z}$ and the common refinement of ${\mathcal A}'_{m,{F_m}}$, for all $m\leq n$. Let $y_n = \sigma({\mathcal V}'_0, \cdots, {\mathcal V}'_n)$, we have the following cases:

\vspace{0.2cm}
   Case 1: If the outcome (i) holds. Then at some stage $n$,  $x_n=s_k$ for some $k$ and $\{s_k\}\in {\mathcal A}_{n}\sim F_n$. Then by definition of the proximal game, all subsequent choices of Player II using $\sigma$ are played inside $X_k\in {\mathcal A}'_{n,{F_n}}$ which, being homeomorphic to a subspace of $\omega_1$, is semi-proximal, so Player I can win the game. 

\vspace{0.2cm}
   Case 2: If the outcome (ii) holds. Then, there is $\alpha\in\omega_1$ such that $a = a_\alpha$ and hence $U_\alpha = a_\alpha\cup\{a_\alpha\}\setminus 2^{\leq n_{F_n}}\in {\mathcal A}_n\sim F_n$. If 
$y_n \in U'_\alpha\cap ((n_{F_n},\omega]\times\alpha_{F_n}+1)$, then the rest of the game will be inside a metric space which is semi-proximal. Otherwise, $y_n\in U'_\alpha\cap ((n_{F_n},\omega])\times(\alpha_{F_n},\omega_1))$ and if this is the case for the rest innings, then $$\bigcap_{n\in\omega}\text{St}(y_n,{ \mathcal V}_n') \subset \bigcap_{n\in\omega}  U'_\alpha\cap ((n_{F_n},\omega])\times(\alpha_{F_n},\omega_1))= \emptyset.$$

\vspace{0.2cm}
   Case 3:  If the outcome (iii) holds. Let $U_n = \text{St}(x_n,{\mathcal A}_n\sim F_n)$.  Then consider $V_n = \text{St}(y_n,{ \mathcal V}_n')$, then $V_n\subset U'_n\cap U$, where $U\in {\mathcal U}_{n_{F_n},F_n}$. We claim that $\bigcap_{n\in\omega} V_n = \emptyset$. To see that, suppose that there is an element $x\in \bigcap_{n\in\omega} V_n$, then $x\in U_n'$, for all $n$. If there is $k$ such that $x\in X_k$, then $X_k\subset U^\prime_n$, for all $n$, and henc $s_k\in U_n$, for all $n$, contradiction. If there is $\alpha$ such that $x = \langle\omega,\alpha+1\rangle$, then $a_\alpha\in U_n$, for all $n$, contradiction.

\vspace{0.2cm}
   Thus $\sigma$ is defeated and hence $X_{A_Z}$ is semi-proximal.
\end{proof}

However, we can show that for finite products of subspaces of $\omega_1$, normality is equivalent to semi-proximality. The case of the product of two subspaces was proved in \cite{Pa}, where the following was shown.

 \begin{lemma} \label{PR} If $A,B\subset\omega_1$, then the following conditions are equivalent:
 \begin{enumerate} 
 \item $A\times B$ is normal.
 \item Either $A$ or $B$ is not stationary, or $A\cap B$ is stationary.
 \item $A\times B$ is semi-proximal.
 \end{enumerate} 
 \end{lemma} 
 Which we extend for finite product of subspaces of $\omega_1$.
 \begin{theorem} \label{ProdFin} The product of finitely many subspaces of $\omega_1$ is semi-proximal if and only if it is normal. 
\end{theorem}

To prove the theorem, we need some preliminary results. 

\begin{lemma}\label{SS} If $X = \prod_{k< n}A_k$ is a finite product of non-empty subspaces of $\omega_1$ such that $X$ is semi-proximal and $\bigcap_{k<n}A_k$ is not stationary,  then there is $k< n$ such that $A_k$ is not stationary.
\end{lemma}

\begin{proof}  We will establish the proof by induction on $n$. The cases $n = 0,1$ are straightforward, and the case where $n = 2$ is given by Theorem \ref{PR}. Assume the statement holds true for $n$ and suppose that $\{A_k:k\leq n\}$ is such that its product $X = \prod_{k\leq n}A_k$ is semi-proximal and $\bigcap_{k\leq n}A_k$ is not stationary.  For each $i\leq n$, define 
$$
X_i = \{\langle\alpha_0,\cdots,\alpha_n\rangle\in X: \exists \beta \,\,\text{such that}\,\,\alpha_k = \beta, \forall k\not=i \}.
$$
For each $i$, $X_i$ is a closed subspace of $X$ and is therefore semi-proximal. Note that each $ X_i$ is homeomorphic to the product of the two sets,  $\left(\bigcap_{k\not=i}A_k\right)\times A_i$. 
Now suppose that $A_k$ is stationary for all $k\leq n$. According to  Theorem \ref{PR}, $\bigcap_{k\not=i}A_k$ is not stationary for all $i\leq n$. By our inductive hypothesis, this implies that $\prod_{k\not=i}A_k$ is not semi-proximal, for all $i$. Consequently, $X_i$ is not semi-proximal which contradicts the assumption that $X$ is semi-proximal. Therefore, there must exist $k\leq n$ for which $A_k$ is not stationary. 
\end{proof} 

\vspace{0.7cm}
   In \cite{PS} Przymusi\'nski proved the following theorem characterizing normality of products with a metric factor. 
\begin{theorem}\label{PS} Let ${\mathcal B}$ be a base for a metrizable space $M$. The product space $X\times M$ is normal if and only if $X$ is normal and for every family $\{F_B:B\in{\mathcal B}\}$ of closed subsets of $X$ such that if $B\subset B^\prime$ then $F_B\subset F_{B^\prime}$ and  for all $z\in M$, $\bigcap\{F_B: z\in B\} = \emptyset$, there exists a family $\{U_B:B\in{\mathcal B}$ of open subsets of $X$ such that $F_B\subset U_B$ and for all $z\in M$, $\bigcap\{U_B:z\in B\} = \emptyset$.
\end{theorem} 

    We now show that it follows from Prymusi\'nski's theorem that if $X$ is normal and countably paracompact, then its product with any countable metrizable space is normal. 

\begin{corollary} \label{NCP} The product of a countably paracompact normal space with a countable metric space is normal.
 \end{corollary} 
\begin{proof}  Let $X$ be a normal and countably paracompact space and let $Y$ be a countable metric space. It is straightforward to show that any countable metrizable space $Y$ has a base ${\mathcal B} = \bigcup_{y\in Y} {\mathcal B}_y$ such that, for all $y\in Y$, ${\mathcal B}_y$ is decreasing local neighborhood base at $y$, and ${\mathcal B}_y\cap {\mathcal B}_{z} = \emptyset$, for $z\neq y$ in $Y$. 

\vspace{0.2cm}
   So we fix such a base and let $\{F_B:B\in{\mathcal B}\}$ be a family of closed subsets of $X$ which satisfies the condition in theorem \ref{PS}. Note that by the monotonicity property of the family of $F_B$'s we have that for every $y\in Y$, $\{F_B: B\in{\mathcal B}_y\}$ is a decreasing sequence of closed subsets of $X$ and since $B_y$ is a local base at $y$ we have $\bigcap\{F_B: B\in{\mathcal B}_y\}\subseteq \bigcap\{F_B: B\in{\mathcal B},\ y\in B\}=\emptyset$. 
Since $X$ is countably paracompact, there exists a family of open sets  $\{U_B:B\in{\mathcal B}_y\}$ in $X$ such that $\bigcap\{U_B:B\in{\mathcal B}_y\} = \emptyset$ and  $F_B\subset U_B$, for all $B\in{\mathcal B}_y$. Since ${\mathcal B}_y$ and ${\mathcal B}_z$ are disjoint for all $y\not=z$, the family $\{U_B:B\in {\mathcal B}\}$ is well defined. To show that $X\times Y$ is normal, it remains to prove that $\{U_B:B\in{\mathcal B}\}$ satisfies the conclusion of theorem \ref{PS}. To see that, let $y\in Y$ then $$\bigcap\{U_B:y\in B\} \subset\bigcap \{U_B:  B\in{\mathcal B}_y\} = \emptyset$$
\end{proof}

\begin{proof}[Proof of Theorem \ref{ProdFin}]  We have already shown in Theorem \ref{NormSemi} that any normal subspace of $\omega_1^n$ is semi-proximal, so we only need to show sufficiency.

   We will prove by induction on $n$ that for any family of $n$ many subspaces of $\omega_1$, if the product is semi-proximal, then the product is normal. The base case $n = 2$, is given by Theorem \ref{PR}. Now assume it is true for $n$ and suppose that $\{A_k:k\leq n\}$ is such that its product $X = \prod_{k\leq n}A_k$ is semi-proximal. Then we have two cases: 

\vspace{0.2cm}
   Case 1: If there exists $k_0<n$ such that $A_{k_0}$ is not stationary, then there exists a club $C$ such that $A_{k_0} = \bigoplus_{\gamma\in C}(\gamma,\gamma^+]\cap A_{k_0}$, where $\gamma^+ = \min C\setminus\gamma$. Then $X =  \bigoplus_{\gamma\in C} X_\gamma$, where 
$$
X_\gamma = \big( (\gamma,\gamma^+]\cap A_{k_0}\big)\times \prod_{k\leq n,k\neq k_0} A_k
$$
Since $\prod_{k< n,k\neq k_0} A_k$ is semi-proximal, then it is normal by induction and hence it is countably paracompact by Theorem \ref{Fl} since shrinkable implies countably paracompact. Therefore, $X_\gamma$ is normal by Corollary \ref{NCP}, for all $\gamma\in C$. Thus, $X$ is normal. 

\vspace{0.2cm}
   Case 2: If $A_k$ is stationary for all $k\leq n$, then $\Delta(X) = \bigcap_{k\leq n}A_k$ is stationary as stated in Lemma \ref{SS}.  Let $H$ and $K$ be two disjoint closed subsets of $X$. Given that $\Delta(X)$ is stationary, there exists an $\alpha\in \Delta(X)$ such that either $H$ or $K$ is a subset of $X^\prime = \bigoplus_{k\leq n} Y_k$, where
$$
Y_k = \alpha+1 \times \prod_{j<k}(A_j\setminus\alpha+1) \times \prod_{k<m\leq n}A_m.
$$ 
To demonstrate this, let’s assume the contrary. Let ${\mathcal M}$ be a countable elementary submodel such that $A_k\in {\mathcal M}$ for all $k\leq n$. We define $\alpha = {\mathcal M}\cap\omega_1\in\Delta(X)$. Consider $(\alpha_i:i\in\omega)$, an increasing sequence that converges to $\alpha$. Then there exists $(h_i:i\in\omega)\subset H$ and $(k_i:i\in\omega)\subset K$ such that $h_i<k_i<\langle\alpha_i,\cdots, \alpha_i\rangle<h_{i+1}$ and hence $(h_i:i\in\omega)$ and $(k_i:i\in\omega)$ would converge to $\langle\alpha,\cdots,\alpha\rangle\in H\cap K$, which contradicts our initial assumption.

\vspace{0.2cm}
   So, suppose $K\subset X^\prime$. Since $ \prod_{j<k}(A_j\setminus\alpha+1) \times \prod_{k<m\leq n}A_m$ is a closed subspace of $X$, it is semi-proximal. By induction,  it is normal, and according to Theorem \ref{Fl}, it is countably paracompact since being shrinkable implies countably paracompact. Therefore,  by Corollary \ref{NCP}, $Y_k$ is normal for all $k\leq n$. Therefore $X^\prime$ is normal, then there exists two disjoint open sets $U^\prime$ and $V$ in $X^\prime$ with $(H\cap X^\prime)\subseteq U^\prime$ and $K \subseteq V$. Since $X^\prime$ is clopen in $X$, then $U^\prime$ and $V$ are open in $X$. Let $U = U^\prime \cup \prod_{k\leq n} (A_k\setminus (\alpha+1))$, then $U$ and $V$ are disjoint open subsets in $X$ which separate $H$ and $K$. Thus, $X$ is normal.
\end{proof}

  The following characterization of countable paracompactness must be known but it seems it has not been explicitly stated before. It follows easily from Corollary \ref{NCP}.

\begin{theorem} A normal space is countably paracompact if and only if its product with any countable metrizable space is normal. 
\end{theorem}

\end{document}